\input amstex
\documentstyle{amsppt}
\magnification 1100
\NoRunningHeads
\NoBlackBoxes
\document

\def\tro{\tilde\rho}
\def\tpi{\tilde\pi}

\def\h{\frak h}

\def\ell{{\text{ell}}}

\def\tA{\tilde A}

\def\R{\Bbb R}
\def\h1{\hat{\bold 1}}

\def\tG{\tilde G}

\def\Ua{U_q(\tilde\g)}
\def\U2{{\Ua}_2}
\def\g{\frak g}

\def\Z{\Bbb Z}
\def\C{\Bbb C}

\def\<{\langle}
\def\>{\rangle}
\def\o{\otimes}
\def\e{\varepsilon}

\topmatter
\title A method of construction of finite-dimensional 
triangular semisimple Hopf algebras\endtitle
\author {\rm Pavel Etingof and Shlomo Gelaki \linebreak
Department of Mathematics\linebreak
Harvard University\linebreak
Cambridge, MA 02138, USA\linebreak
e-mail: etingof\@math.harvard.edu\linebreak shlomi\@math.harvard.edu}
\endauthor

\endtopmatter
\centerline{June 12, 1998}
\head {\bf Introduction}\endhead
The goal of this paper is to 
give a new method of constructing finite-dimensional semisimple triangular 
Hopf algebras, including minimal ones which are non-trivial (i.e. not group 
algebras). The paper shows that such Hopf algebras are quite abundant. 
It also discovers an unexpected connection of such Hopf algebras with 
bijective 1-cocycles on finite groups and set-theoretical solutions 
of the quantum Yang-Baxter equation defined by Drinfeld \cite{Dr1}.

Finite-dimensional triangular Hopf algebras were studied by several authors 
(see e.g. [CWZ,EG,G,M]). In [EG] the authors prove that any 
finite-dimensional semisimple 
triangular Hopf algebra over an algebraically closed field of 
characteristic $0$ (say $\C$) is obtained from a group algebra after 
twisting 
its comultiplication in the sense of Drinfeld [Dr2]. Twists are easy to 
construct for abelian groups $A,$ they are just 2-cocycles for $A^*$ with 
values in $\C^*.$ A general simple construction of triangular semisimple 
Hopf algebras which are non-trivial, is 
the following: take a non-abelian group $G,$ an abelian subgroup of it $A,$
and a twist $J\in \C[A]\o \C[A]$ which does not commute with $g\o g$ for 
all $g\in G,$ and twist $\C[G]$ by $J$ to obtain $\C[G]^J.$
Examples of such Hopf algebras were constructed by the second author in [G]. 
They are Hopf algebras of dimension $pq^2$ where $p$ and $q$ are any 
prime numbers so that $q$ divides $p-1.$ It was also proved in [G] that 
the dual of the Drinfeld 
double of these Hopf algebras is triangular. Nevertheless Hopf 
algebras which are constructed in this way are not minimal, and their 
minimal Hopf subalgebras are trivial. In fact, as far as we know, in the 
literature there are no non-trivial semisimple minimal
triangular Hopf algebras. A natural question thus arose: Are 
there finite-dimensional non-trivial minimal semisimple triangular Hopf 
algebras? In this paper we describe a method for constructing such Hopf 
algebras.

The paper is organized as follows. First, we show how to construct twists for 
certain solvable non-abelian groups by iterating twists of their 
abelian subgroups, and thus obtain non-trivial semisimple triangular Hopf 
algebras. Second, we show that in some cases this construction gives
non-trivial semisimple {\it minimal} triangular Hopf algebras. Finally, we 
show how any non-abelian group which admits a bijective 1-cocycle gives 
rise to a non-trivial semisimple {\it minimal} triangular Hopf algebra. 
Such non-abelian groups exist in abundance and were constructed in [ESS] 
in connection with set-theoretical solutions to the quantum Yang-Baxter 
equation.  

We shall work over the field of complex numbers $\C,$ although it can be 
replaced by any algebraically closed field of characteristic $0.$  
\head {\bf 1. Twists}\endhead

Recall Drinfeld's notion of a twist for Hopf algebras. 

\proclaim{Definition 1.1 \cite{Dr2}}
Let ${\Cal A}$ be a Hopf algebra over a field $k.$ A twist for ${\Cal A}$ 
is an invertible element $J\in {\Cal A}\o {\Cal A}$ which satisfies 
$$
(\Delta\o I)(J)J_{12}=(I\o \Delta)(J)J_{23}\ and \ (\e\o I)(J)=(I\o \e)(J)=1
\tag 1.1
$$
where $I:{\Cal A}\to {\Cal A}$ is the identity map.
\endproclaim

Given a twist $J$ for ${\Cal A}$, we can construct a new 
Hopf algebra ${\Cal A}^J$, which is the same as ${\Cal A}$ as an algebra, 
with coproduct $\Delta^J$ given by
$$
\Delta^J(x)=J^{-1}\Delta(x)J,\ x\in {\Cal A}.\tag 1.2
$$
If $({\Cal A},R)$ is quasitriangular then so is ${\Cal A}^J$ with the 
R-matrix $$
R^J=J_{21}^{-1}RJ,\tag 1.3
$$
and if $({\Cal A},R)$ is triangular then so is $({\Cal A}^J,R^J).$
Note that twists can be composed; that is, if $J$ is a twist for ${\Cal A}$ 
and $J'$ is a twist 
for ${\Cal A}^J$ then $JJ'$ is a twist for ${\Cal A}$, and 
${\Cal A}^{JJ'}=({\Cal A}^J)^{J'}$. 

Now let $({\Cal A},R)$ be a finite-dimensional semisimple triangular Hopf 
algebra 
over $\C$. Assume that the Drinfeld element $u$ equals 1 (this can always 
be attained by a simple modification of $R$, see \cite{EG}). It was shown in 
[EG, Theorem 2.1] that in this case there exists a finite group $G$ such that 
${\Cal A}=\C[G]^J$ as a triangular Hopf algebra, where $J$ is a suitable 
twist. Thus, construction of triangular semisimple finite-dimensional 
Hopf algebras reduces to construction of twists for group algebras. 

\proclaim{Remark 1.2} {\rm In \cite{CWZ}, the authors generalize a 
number of results on commutative algebras to the case of quantum commutative 
algebras in the tensor category of $\Cal A$-modules, where $\Cal A$ is a 
finite-dimensional semisimple triangular Hopf algebra. Let us point out 
that [EG, Theorem 2.1] is very useful for proving such generalizations. For 
example, the main result of \cite{CWZ} is 
Theorem 4.7, which in particular states that if the 
characteristic of the ground field is $0$ or $>dim{\Cal A},$ then 
any quantum commutative algebra $R$ in $\Cal A$-mod is integral over 
its subalgebra of invariants $R^{\Cal A}$. Let us show 
how [EG, Theorem 2.1] enables to give an easy proof of this result in 
characteristic $0,$ when the Drinfeld element $u$ of $\Cal A$ equals $1.$ 

Let $\alpha\in R$ and pick a finite-dimensional $\Cal A$-submodule $V$ of 
$R$ so that 
$\alpha\in V.$ We wish to show that the images of $K_i=\oplus_{j=0}^iR^{\Cal 
A}\o V^{\o j}$ under multiplication $K_i\to R$ stabilize. This will imply 
integrality of 
$R$ over $R^{\Cal A}$, as then the 
 algebra $R^{\Cal A}[V]$ is a faithful module 
over $R^{\Cal A}[\alpha]$ which is finitely generated
over $R^{\Cal A}$ (see [L, p. 356]).

Theorem 2.1 in [EG] states that there exists a fiber functor $F:{\Cal A} 
-\text{mod}\to Rep(G)$ which is an equivalence 
of $\Cal A$-mod to the category of representations of some finite group $G$. 
Thus, $R_0=F(R)$ is an ordinary commutative algebra with a $G$-action, 
and $F(R^{\Cal A})=R_0^G.$
Now, we know that $R_0$ is integral over $R_0^G$ of degree $|G|$
(since any $x\in R_0$ satisfies $P(x)=0$ where $P(y)=\prod_{g\in G}(y-gx)$). 
This implies that if 
$U$ is any finite-dimensional subobject (i.e. a $G$-submodule) of $R_0$, 
then the images of $K_i=\oplus_{j=0}^iR_0^G\o U^{\o j}$ 
under the multiplication map $K_i\to R_0,$ stabilize (after $i=|G|$). 
In particular, it is true for $U=F(V)$. But,
this property is categorical, i.e. preserved by tensor functors. 
Therefore, it remains valid if we replace $R_0$ with $R$, $R_0^G$ with 
$R^{\Cal A}$, and $U$ with $V$, which implies the claim. $\square$ }
\endproclaim

The following definition is due to Radford [R].

\proclaim{Definition 1.3} Let $({\Cal A},R)$ be a quasitriangular Hopf 
algebra, and write $R=\sum_i x_i\o y_i$ in the shortest possible way. 
Then $({\Cal A},R)$ is called minimal if ${\Cal A}$ is 
generated by its Hopf subalgebras $sp\{x_i\}$ and $sp\{y_i\}.$ 
\endproclaim

It is easy to check that if $({\Cal A},R)$ is triangular, then it is 
minimal if and only if $R$ defines a nondegenerate bilinear form on 
${\Cal A}^*$ (see e.g. [G]). A twist $J$ for $\C[G]$ for which 
$(\C[G]^J,J_{21}^{-1}J)$ is minimal triangular is said to be {\it 
minimal}. By \cite{EG}, for any finite-dimensional semisimple 
triangular Hopf algebra $({\Cal A},R)$ there exist finite groups 
$H\subset G$ and  a minimal twist $J$ for $\C[H]$ such that 
${\Cal A}=\C[G]^J$. Therefore, construction of finite-dimensional 
semisimple triangular Hopf algebras reduces to construction 
of minimal ones.

For a general group $G$, minimal twists for $\C[G]$ are quite difficult 
to construct. However, if $G$ is abelian, they can be constructed easily 
(if they exist). Namely, a twist for $\C[G]$, when regarded as a 
$\C^*$-valued function on $G^*\times G^*$, 
is the same thing as a 2-cocycle of the group $G^*$ with 
coefficients in $\C^*$. In many cases, such a twist is minimal. 

Our goal here is to use this simple construction to 
construct twists in the case of a non-abelian $G$. 

\head {\bf 2. Construction of twists for non-abelian groups}\endhead

In this section we give a method of constructing twists 
for non-abelian groups.
Let $G,A$ be groups and $\rho :G\to Aut(A)$ a homomorphism. 
For brevity we will write $\rho(g)(x)=gx$ for $g\in G,x\in A.$ 

\proclaim{Definition 2.1}
By a 1-cocycle of $G$ with coefficients in $A$ we mean a map 
$\pi:G\to A$ which satisfies the equation 
$$
\pi(gg')=\pi(g)(g\pi(g')),\ g,g'\in G.\tag 2.1
$$ 
\endproclaim

We will be interested in the case when $\pi$ is {\it a bijection}, 
because of the following proposition. 

\proclaim{Proposition 2.2} Let $G,A$ be groups, $\pi:G\to A$ a 
bijective 1-cocycle, and $J$ a twist for $\C[A]$ which is $G$-invariant. 
Then $\bar J:=(\pi^{-1}\o \pi^{-1})(J)$ 
satisfies (1.1). Thus, if $\bar J$ is invertible then it is a twist for 
$\C[G]$. \endproclaim

\demo{Proof} It is obvious that the second equation of (1.1) is satisfied 
for $\bar J$. So we only have to prove the  
first equation of (1.1) for $\bar J$. Let $J=\sum a_{xy}x\o y$. Then 
$$
\gather
(\pi\o\pi\o\pi)((\Delta\o I)(\bar J)\bar J_{12})=\\
\sum_{x,y,z,t\in A}a_{xy}a_{zt}\pi(\pi^{-1}(x)\pi^{-1}(z))\o
\pi(\pi^{-1}(x)\pi^{-1}(t))\o \pi(\pi^{-1}(y))=\\
\sum_{x,y,z,t\in A}a_{xy}a_{zt}x(\pi^{-1}(x)z)\o
x(\pi^{-1}(x)(t))\o y.\tag 2.2\endgather
$$
Using the $G$-invariance of $J$, we can remove the $\pi^{-1}(x)$ 
in the last expression and get 
$$
(\pi\o\pi\o\pi)((\Delta\o I)(\bar J)\bar J_{12})=
(\Delta\o I)(J)J_{12}.\tag 2.3
$$
Similarly, 
$$
(\pi\o\pi\o\pi)((I\o \Delta)(\bar J)\bar J_{23})=
(I\o \Delta)(J)J_{23}.\tag 2.4
$$
But $J$ is a twist, so the right hand sides of (2.3) and (2.4) are equal. 
Since $\pi$ is bijective, this implies equation (1.1) for $\bar J$. 
$\square$\enddemo
 
Let $G_1,...,G_n=G$, $H_1,...,H_n$ be finite groups,
and $\rho_{i-1}:G_{i-1}\to Aut(H_i)$, $i=2,...,n$, be 
homomorphisms, such that $G_1=H_1$, and $G_i=G_{i-1}\ltimes H_i$, 
where the semidirect product is made using the action $\rho_{i-1}$. 
We will call such data a {\it hierarchy of length $n$}. 
Note that $G$ canonically contains $G_1,...,G_n$ and $H_1,...,H_n$ 
as subgroups. Let $A=H_1\times\cdots\times H_n$, and 
$\rho:G\to Aut(A)$ be defined by $\rho(x_n\cdots x_1)(y_1,...,y_n)=
(y_1,\rho_1(x_1)y_2,...,\rho_{n-1}(x_{n-1}\cdots x_1)y_n).$

\proclaim{Proposition 2.3}
The map $\pi:G\to A$ given by $\pi(x_n\cdots x_1)=(x_1,...,x_n)$, 
$x\in H_i$ is a bijective 1-cocycle. 
\endproclaim  

\demo{Proof} Straightforward. 
$\square$\enddemo

In the situation of Proposition 2.3 suppose further that
$J_i$ is a twist for $\C[H_i]$, $i=1,...,n$,
and that $J_i$ is invariant under $\rho_{i-1}(G_{i-1})$, $i=2,...,n.$
We will call such data a {\it twist hierarchy of length $n$}.
In this situation, the element $J=J_1\o\cdots\o J_n$ is a twist for $\C[A]$. 
Therefore by
Proposition 2.2, $\bar J=(\pi^{-1}\o\pi^{-1})(J)$ satisfies (1.1). 

\proclaim{Proposition 2.4}
$\bar J=J_n\cdots J_1,$ and is a twist for $\C[G].$
\endproclaim

\demo{Proof} It is clear that $\bar J$ is invertible and ${\bar 
J}^{-1}=J_1^{-1}\cdots J_n^{-1}.$ Therefore, $\bar J$ is a twist for $\C[G].$
$\square$\enddemo

\head {\bf 3. Minimal triangular Hopf algebras arising from symplectic 
hierarchies}\endhead

In this section we show that the twist constructed in Proposition 2.4 is 
often minimal.

Let $H$ be a finite abelian group, and
$(,):H\times H^*\to 2\pi i(\R/\Z)$ be the standard pairing. 
By a {\it symplectic structure on $H$} we mean an isomorphism
$B:H\to H^*$ such that $B^*=-B$. If $B$ is a symplectic structure
then the bilinear form $(x,By)$ on $H$ will be denoted by $\<x,y\>$. 
A symplectic structure defines an element
$$
J_B=|H|^{-1}\sum_{x,y\in H}e^{\<x,y\>}x\o y\in \C[H]\o \C[H].\tag 3.1
$$
This element is invertible with $J_B^{-1}=|H|^{-1}\sum_{x,y\in 
H}e^{-\<x,y\>}x\o y,$ so it is a twist for $\C[H].$

It is clear that an automorphism of $H$ fixes $J_B$ if and only if it fixes 
$B$. 

If $H$ is an abelian group of odd order, then taking the square
is an isomorphism $H\to H$. 
For $x\in H$, we will denote the preimage of $x$ under this isomorphism by 
$x^{1/2}$. 

If $|H|$ is odd then it is straightforward to check that the R-matrix of 
$\C[H]^{J_B}$ is 
$R_B=(J_B)_{21}^{-1}J_B=|H|^{-1}\sum e^{\<x,y^{1/2}\>}x\o y.$ The 
matrix $c_{xy}=|H|^{-1}e^{\<x,y^{1/2}\>}$ is nondegenerate; its inverse 
is $e^{-\<x,y^{1/2}\>}$ by the inversion formula for Fourier transform. 
So if $|H|$ is odd then $J_B$ is minimal.

Consider now a twist hierarchy of length $n,$ of the following 
form: $H_i$ are abelian of odd order, and $J_i=J_{B_i}$, 
where $B_i$ is a symplectic structures on $H_i$ for all $i.$ Such a hierarchy 
will be called {\it symplectic}. 

\proclaim{Theorem 3.1} For a symplectic hierarchy, the twist 
$\bar J=J_n\cdots J_1\in \C[G]\o \C[G]$ is minimal, i.e. the
triangular Hopf algebra $(\C[G]^{\bar J}, {\bar J}_{21}^{-1}{\bar J})$ is 
minimal. Its universal R-matrix has the form 
$$
R=|G|^{-1}\sum_{x_i,y_i\in H_i}
e^{\sum_{j=1}^n\<\rho_{j-1}(x_{j-1}^{1/2}\cdots x_1^{1/2})^{-1}x_j,
\rho_{j-1}(y_{j-1}^{1/2}\cdots y_1^{1/2})^{-1}y_j^{1/2}\>}
x_n\cdots x_1\o y_n\cdots y_1.\tag 3.2
$$
\endproclaim

Before we prove Theorem 3.1 we need the following lemma.

\proclaim{Lemma 3.2}
Let $K$ be an abelian group of odd order 
with a symplectic form, acting symplectically
on another abelian group $H$ of odd order 
with a symplectic form via $\rho:K\to Aut(H).$ Let
$$
S_{K,H}(x,y)=\sum_{z,z',t,t'\in K}e^{\<z,t\>+\<z',t'\>+\<\rho(z^{-1})x, 
\rho(t^{-1})y^{1/2}\>}zz'\o tt'.\tag 3.3
$$
Then,
$$
S_{K,H}(x,y)=|K|\sum_{u,v\in K}e^{\<u,v^{1/2}\>+
\<\rho(u^{-1/2})x,\rho(v^{-1/2})y^{1/2}\>}u\o v.\tag 3.4
$$
\endproclaim

\demo{Proof}
Set in (3.3) $u=zz',v=tt'$. Then (3.3) transforms to the 
form 
$$
S_{K,H}(x,y)=\sum_{z,u,t,v\in K}e^{2\<z,t\>+\<u,v\>-\<u,t\>-\<z,v\>+
\<\rho(z^{-1})x, 
\rho(t^{-1})y^{1/2}\>}u\o v.\tag 3.5
$$
Introducing a new variable $b=zt^{-1},$ we can sum over $b$ instead of 
$z$ in (3.5). This yields
$$
S_{K,H}(x,y)=\sum_{b,u,t,v\in K}e^{2\<bt,t\>+\<u,v\>-\<u,t\>-\<bt,v\>+
\<\rho(t^{-1})\rho(b^{-1})x, 
\rho(t^{-1})y^{1/2}\>}u\o v.\tag 3.6
$$
Using the skew symmetry of $\<,\>$ (which implies $\<t,t\>=0$ because 
of odd order) and the invariance of the symplectic form on $H$ under $K$, we 
get 
$$
S_{K,H}(x,y)=\sum_{b,u,t,v\in K}e^{\<u,v\>+\<b^2u^{-1}v,t\>-\<b,v\>+
\<\rho(b^{-1})x,y^{1/2}\>}u\o v.\tag 3.7
$$
Now we can sum over $t$, using the identity 
$\sum_t e^{\<a,t\>}=|K|\delta_{a1}$. This yields
$$
S_{K,H}(x,y)=|K|\sum_{u,v\in K}e^{\<u^{1/2},v\>+
\<\rho(u^{1/2}v^{-1/2})x,y^{1/2}\>}u\o v,\tag 3.8
$$
and we are done. $\square$\enddemo

\noindent {\it Proof of Theorem 3.1.}
The first statement clearly follows from the second one. 
Indeed, let $R=\sum_{g,h\in G}a_{gh}g\o h$. It follows from (3.2) 
that by permutation of rows and columns the matrix $a=(a_{gh})$ can be 
transformed into $b=(b_{gh})$, where 
$$
b_{x_n\cdots x_1,y_n\cdots y_1}=|G|^{-1}e^{\sum\<x_i,y_i\>}.\tag 3.9
$$
This matrix is a tensor product $b^n\o\cdots \o b^1$, where 
$b^i_{x,y}=|H_i|^{-1}e^{\<x,y\>}$, $x,y\in H_i$. 
As we mentioned, 
the matrices $b_i$ are nondegenerate, therefore so are $b$ and $a$. 

We now calculate $R$. It is straightforward to check 
that the twist $J_i$ satisfies 
$J_i^{21}=J_i^{-1}$ for all $i.$ Therefore, $R=J_1\cdots 
J_{n-1}J_n^2J_{n-1}\cdots J_1$. It is also straightforward to check that 
$J_B^2=J_{B'}$ for any $B,$ where $B'x=Bx^{1/2}$. Therefore, we have: 
$$
\gather
R=\frac{|H_n|}{|G|^2}\sum_{z_i,t_i,x_n',y_n',z_i',t_i'}
e^{\sum_i\<z_i,t_i\>+\sum_i\<z_i',t_i'\>+\<x_n',(y_n')^{1/2}\>}  
\times\\ 
z_1\cdots z_{n-1}x_n'z_{n-1}'\cdots z_1'\o 
t_1\cdots t_{n-1}y_n't_{n-1}'\cdots t_1' \tag 3.10\endgather $$ 
Moving $x_n'$ and $y_n'$ to the left, we get from (3.10):
$$
\gather
R=\frac{|H_n|}{|G|^2}\sum_{z_i,t_i,x_n',y_n',z_i',t_i'}
e^{\sum_i\<z_i,t_i\>+\sum_i\<z_i',t_i'\>+\<x_n',(y_n')^{1/2}\>}\times\\
x_nz_1\cdots z_{n-1}z_{n-1}'\cdots z_1'\o 
y_nt_1\cdots t_{n-1}t_{n-1}'\cdots t_1'\tag 3.11 \endgather $$ 
where $x_n=\rho_{n-1}(z_1\cdots z_{n-1})x_n'$ and 
$y_n=\rho_{n-1}(t_1\cdots t_{n-1})y_n'$. Now we can replace summation over 
$x_n',y_n'$ with summation over $x_n,y_n$ and get
$$
\gather
R=\frac{|H_n|}{|G|^2}\sum_{z_i,t_i,x_n,y_n,z_i',t_i'}
e^{\sum_i\<z_i,t_i\>+\sum_i\<z_i',t_i'\>+\<\rho_{n-1}(z_{n-1}^{-1}\cdots 
z_1^{-1})x_n, \rho_{n-1}(t_{n-1}^{-1}\cdots t_1^{-1})y_n^{1/2}\>}\times\\
x_nz_1\cdots z_{n-1}z_{n-1}'\cdots z_1'\o 
y_nt_1\cdots t_{n-1}t_{n-1}'\cdots t_1'.\tag 3.12 \endgather
$$ 
Thus, the statement of the theorem is equivalent 
to the identity
$$
\gather
\frac{|H_n|}{|G|}\sum_{z_i,t_i,z_i',t_i'}
e^{\sum_i\<z_i,t_i\>+\sum_i\<z_i',t_i'\>+\<\rho_{n-1}(z_{n-1}^{-1}\cdots 
z_1^{-1})x_n, \rho_{n-1}(t_{n-1}^{-1}\cdots t_1^{-1})y_n^{1/2}\>}\times\\
z_1\cdots z_{n-1}z_{n-1}'\cdots z_1'\o t_1\cdots t_{n-1}t_{n-1}'\cdots t_1'=\\
\sum_{x_i,y_i\in H_i,1\le i\le n-1}
e^{\sum_{j=1}^n\<\rho_{j-1}(x_{j-1}^{1/2}\cdots x_1^{1/2})^{-1}x_j,
\rho_{j-1}(y_{j-1}^{1/2}\cdots y_1^{1/2})^{-1}y_j^{1/2}\>}\times\\
x_{n-1}\cdots x_1\o y_{n-1}\cdots y_1.\tag 3.13 
\endgather
$$ 
So it remains to prove (3.13).

Denote the left hand side of (3.12) by $T_{H_1,...,H_n}(x_n,y_n)$.
Using Lemma 3.2 for $K=H_{n-1}$, $H=H_n$, 
$x=x_n$, $y=y_n$, we get
$$
\gather
T_{H_1,...,H_n}(x_n,y_n)=
\frac{|H_n||H_{n-1}|}{|G|}\sum_{z_i,t_i,z_i',t_i',x_{n-1}',y_{n-1}'}
e^{\sum_i\<z_i,t_i\>+\sum_i\<z_i',t_i'\>}\times \\ 
e^{\<x_{n-1}',(y_{n-1}')^{1/2}\>
+\<\rho_{n-1}((x_{n-1}')^{-1/2})\rho_{n-1}(z_{n-2}^{-1}\cdots z_1^{-1})x_n,
\rho_{n-1}((y_{n-1}')^{-1/2})\rho_{n-1}(t_{n-2}^{-1}\cdots 
t_1^{-1})y_n^{1/2}\>} \times\\
z_1\cdots z_{n-2}x_{n-1}'z_{n-2}'\cdots z_1'\o t_1\cdots t_{n-2}
y_{n-1}'t_{n-2}'\cdots t_1'.\tag 3.14
\endgather
$$ 
(here $i$ runs from $1$ to $n-2$). 
Let us now move $x_{n-1}'$ and $y_{n-1}'$ to the left. Then we get
$$
\gather
T_{H_1,...,H_n}(x_n,y_n)= \frac{|H_n||H_{n-1}|}{|G|}\\
\sum_{z_i,t_i,z_i',t_i',x_{n-1},y_{n-1}}
e^{\sum_i\<z_i,t_i\>+\sum_i\<z_i',t_i'\>
+\<\rho_{n-2}(z_{n-2}^{-1}\cdots z_1^{-1})x_{n-1},
\rho_{n-2}(t_{n-2}^{-1}\cdots t_1^{-1})y_{n-1}^{1/2}\>
}\times \\ 
e^{
\<\rho_{n-1}(z_{n-2}^{-1}\cdots z_1^{-1})\rho_{n-1}(x_{n-1}^{-1/2})x_n,
\rho_{n-1}(t_{n-2}^{-1}\cdots t_1^{-1})\rho_{n-1}
(y_{n-1}^{-1/2})y_n^{1/2}\>}\times\\
x_{n-1}z_1\cdots z_{n-2}z_{n-2}'\cdots z_1'\o y_{n-1}t_1\cdots t_{n-2}
t_{n-2}'\cdots t_1'.\tag 3.15
\endgather
$$ 
Thus, 
$$
\gather
T_{H_1,...H_{n-1},H_n}(x_n,y_n)=
\sum_{x_{n-1},y_{n-1}\in H_{n-1}}(x_{n-1}\o y_{n-1})\times \\
T_{H_1,...,H_{n-2},H_{n-1}\times H_n}(x_{n-1}\rho_{n-1}(x_{n-1}^{-1/2})x_n, 
y_{n-1}\rho_{n-1}(y_{n-1}^{-1/2})y_n).\tag 3.16
\endgather
$$

Now we can easily prove (3.13) by induction on $n$. The base of induction 
($n=1$) is obvious, and the induction step follows directly from (3.16). 
This concludes the proof of the theorem.
$\square$

\head {\bf 4. The double of a bijective 1-cocycle and minimal 
triangular Hopf algebras}\endhead

In fact, the method of Proposition 2.2 leads to more examples of triangular 
structures than described in Sections 2,3. Namely, 
given a quadruple $(G,A,\rho,\pi)$ as in Section 2, 
such that $A$ is abelian, define $\tG=G\ltimes A^*$, $\tA=A\times A^*$, 
$\tro:\tG\to Aut(\tA)$ by $\tro(g)=\rho(g)\times\rho^*(g)^{-1}$, 
and $\tpi:\tG\to \tA$ by 
$\tpi(a^*g)=\pi(g)a^*$, $a^*\in A^*$, $g\in G$. It is straightforward to 
check that $\tpi$ is a bijective 1-cocycle. We call the quadruple 
$(\tG,\tA,\tro,\tpi)$ the {\it double of $(G,A,\rho,\pi)$}. 

Consider $J\in \C[\tA]\o \C[\tA]$ given by 
$$
J=|A|^{-1}\sum_{x\in A,y^*\in A^*}e^{(x,y^*)}x\o y^*.\tag 4.1
$$
It is straightforward to check that $J$ is a twist, and that it is 
$G$-invariant. This allows to construct the corresponding element 
$$
\bar J=
|A|^{-1}\sum e^{(x,y^*)}\pi^{-1}(x)\o y^*.\tag 4.2
$$

\proclaim{Proposition 4.1} $\bar J$ is invertible, and 
$$
\bar J^{-1}=|A|^{-1}
\sum_{z\in A,t^*\in A^*}e^{-(z,t^*)}\pi^{-1}(T(z))\o t^*,\tag 4.3
$$
where $T:A\to A$ is a bijective map (not a homomorphism, in general) 
defined by $\pi^{-1}(x^{-1})\pi^{-1}(T(x))=1$. 
\endproclaim

\demo{Proof} Denote the right hand side of (4.3) by $J'$. 
We need to check that $J'=\bar J^{-1}$. It is enough to check it 
after evaluating any $\alpha\in A$ on the second component of both sides. 
We have 
$$
\gather
(1\o\alpha)({\bar J})=|A|^{-1}\sum_{x,y^*}e^{(x\alpha,y^*)}\pi^{-1}(x)=
\pi^{-1}(\alpha^{-1}),\\
(1\o\alpha)(J')=|A|^{-1}\sum_{x,y^*}e^{-(x\alpha^{-1},y^*)}\pi^{-1}(T(x))=
\pi^{-1}(T(\alpha)).\tag 4.4\endgather
$$
This concludes the proof of the proposition. 
$\square$\enddemo

We can now prove:

\proclaim{Theorem 4.2} Let $\bar J$ be as in (4.2). Then $\bar J$ is a twist 
for $\C[\tilde G],$ and
it gives rise to a minimal triangular Hopf algebra $\C[\tilde G]^{\bar J}$, 
with universal R-matrix 
$$
R=|A|^{-2}\sum_{x,y\in A,x^*,y^*\in A^*}e^{(x,y^*)-(y,x^*)}
x^*\pi^{-1}(x)\o \pi^{-1}(T(y))y^*.\tag 4.5
$$
\endproclaim

\demo{Proof} The first two statements follow directly from 
Propositions 2.2,4.1. The minimality follows from the fact 
that $\{x^*\pi^{-1}(x)|x^*\in A^*,x\in A\}$ and 
$\{\pi^{-1}(T(y))y^*|y\in A,y^*\in A^*\}$ are bases of 
$\C[\tG],$ and the fact that the matrix $c_{xx^*,yy*}=
e^{(x,y^*)-(y,x^*)}$ is invertible (because it is proportional to
the matrix of Fourier transform on $A\times A^*$).
$\square$\enddemo

Thus, every bijective 1-cocycle $\pi :G\to A$ gives rise to a minimal 
triangular 
structure on $\C[G\ltimes A^*]$. So it remains to construct a supply 
of bijective 1-cocycles. This was done in \cite{ESS}. 
The theory of bijective 1-cocycles was developed 
in \cite{ESS}, because it was found that they correspond to set-theoretical
solutions of the quantum Yang-Baxter equation. In particular, many 
constructions of these 1-cocycles were found. We refer the reader to 
\cite{ESS} for further detail. 

\proclaim{Example 4.3} {\rm Let $G=S_3$ be the permutation group of 
three letters, and $A=\Z_2\times \Z_3.$
Define an action of $G$ on $A$ by $s(a,b)=(a,(-1)^{sign(s)}b)$ for $s\in 
G,$ $a\in \Z_2$ and $b\in \Z_3.$ Define a bijective 1-cocycle 
$\pi=(\pi _1,\pi _2):G\to A$ as follows: $\pi _1(s)=0$ if $s$ is even and 
$\pi _1(s)=1$ if $s$ is odd, and $\pi _2(id)=0,$ $\pi _2((123))=1,$ 
$\pi _2((132))=2,$ $\pi _2((12))=2,$ $\pi _2((13))=0$ and $\pi 
_2((23))=1.$ 
Then by Theorem 4.2, $\C[\tilde G]^{\bar J}$ is a non-trivial minimal 
triangular Hopf algebra of 
dimension $36.$ This is in fact a special case of the construction in 
Proposition 2.4.} 
\endproclaim

\proclaim{Remark 4.4} {\rm Note that in \cite{ESS}, the 1-cocycle equation 
is $\pi(gg')=\pi(g')((g')^{-1}\pi(g))$, rather than (2.1). 
However, these equations are equivalent and related by the 
transformation $g\to g^{-1}$.} 
\endproclaim
\proclaim{Remark 4.5}
{\rm The Lie-theoretical analogue of the theory of bijective 
1-cocycles for finite groups 
is the theory of left invariant affine space structures on Lie groups, 
which has been intensively discussed in the literature
(see \cite{Bu} and references therein). We note that many methods 
of constructing such structures on nilpotent Lie groups 
(which are based on Lie algebra theory and the Campbell-Hausdorff formula) 
work over a finite field and therefore apply to p-groups.} 
\endproclaim
In conclusion let us formulate some questions which are inspired by the 
above results. 
\proclaim{Question 4.6}
{\rm 1. In the situation of Proposition 2.2, is $\bar J$ always 
invertible? If yes and if $J$ is minimal, is $\bar J$ minimal?

2. In the situation of Proposition 2.4, assume that the twists 
$J_1,...,J_n$ 
are minimal. Does this imply that ${\bar J}=J_n\cdots J_1$ is minimal? 

3. If (in the setting of Section 1)
the Hopf algebra $\C[G]^J$ is minimal triangular, does $G$ have to be 
solvable (cf. Remark 4.7.)?} 
\endproclaim
\proclaim{Remark 4.7} {\rm It was shown in \cite{ESS}
that if a bijective 1-cocycle on $G$ exists then $G$ has to be solvable. 
However, this theorem relies on a non-trivial theorem of P. Hall 
in finite group theory, saying that if $G$ with $|G|=\prod p_i^{n_i}$ has 
a subgroup of index $p_i^{n_i}$ for any $i$ then $G$ is solvable.} 
\endproclaim

\end